\documentclass{siamltex}
\usepackage{graphicx}
\usepackage{subfig}
\usepackage{tabularx}
\usepackage{amsmath}
\usepackage{latexsym}
\usepackage{amssymb}
\usepackage{color}
\usepackage{amsfonts}
\usepackage{mathdots}
\usepackage{ulem}
\usepackage{algorithm,algorithmic}
\usepackage{fancyhdr}
\usepackage{multirow}
\usepackage{epsfig}
\usepackage{colortbl}
\usepackage{xcolor}

\makeatletter
\newcommand{\ssymbol}[1]{^{\@fnsymbol{#1}}}
\makeatother

\usepackage[mathscr]{euscript}

\title{Multilinear Discriminant Analysis using a new family of tensor-tensor products}

\author{F.Dufrenois \thanks{ LISIC,50 rue F. Buisson, ULCO Calais, France} \and A. El Ichi\footnotemark [3] \thanks{Laboratoire de Mathématiques, Informatique et Applications, S\'ecurit\'e de l'Information LABMIA-SI, University Mohamed V, Rabat Morocco }  \and    K. Jbilou\footnotemark[1] \thanks{LMPA, 50 rue F. Buisson, ULCO Calais, France; Mohammed VI Polytechnic University, Green City, Morocco; jbilou@univ-littoral.fr }}

\begin{document}

\maketitle

\begin{abstract}    
	Multilinear Discriminant Analysis (MDA) is a powerful dimension reduction method specifically formulated to deal with tensor data. Precisely, the goal of MDA  is to find mode-specific projections that optimally separate tensor data from different classes. However, to solve this task, standard MDA methods use alternating optimization heuristics involving the computation of a succession of tensor-matrix products. Such approaches are most of the time difficult to solve and not natural, highligthing the difficulty to formulate this problem in fully tensor form. In this paper, we propose to solve multilinear discriminant analysis (MDA) by using the concept of transform domain (TD) recently proposed in \cite{Kilmer2011}. We show here that moving MDA to this specific transform domain make its resolution easier and more natural. More precisely, each frontal face of the transformed tensor is processed independently to build a separate optimization sub-problems easier to solve. Next, the obtained solutions are converted into projective tensors by inverse transform. By considering a large number of experiments, we show the effectiveness of our approach with respect to existing MDA methods.
\end{abstract}

\noindent {\bf Keywords.}  Krylov subspaces, Linear tensor equations, Tensor L-product.


\medskip

\noindent {\bf AMS Subject Classification} { 65F10, 65F22. }


\section{introduction}

Linear Discriminant Analysis (LDA) is a supervised dimensionality
reduction tool allowing the classification to multiple categories
in datasets. There are used in numerous areas as diverse as speech
and music classification \cite{Cor2005}, video classification \cite{Pan2005}, outlier
detection \cite{Rot2006}, supervised novelty detection \cite{Dufrenois2016,Bodesheim2013}, etc..
Considering that each data is labeled, the goal of LDA is to find
a set of projections which maximizes the between-class scatter while
minimizing the within-class scatter. In the litterature, this objective
is commonly formulated in two different ways as solving: a \textbf{trace
	ratio problem} which is typically nonconvex and need an iterative
optimization procedure or a \textbf{ratio trace problem} which is
inexact but equivalent to a generalized eigenvalue problem (GEP). 

When dealing with high order data such as grayscale images, RGB images,
multispectral images, ... a conventional practice is to vectorize
the whole data set before applying LDA. This preprocessing involves
systematically an increase of the dimensionality of the data sample
and may result in singularity problems commonly referred as the small
sample size problem (SSS). 

A way to solve this question is to adopt the tensor representation
which allows to preserve the natural multidimensional form of the
data while reducing implicitly the dimensionality of the data. A n-order
tensor can be seen as a hyper-parallepiped with n sides and where
each side represents a "mode". A tensor generalizes thus the notions
of matrix (2-order tensor) and vector (1-order tensor). 

In this framework, to perform LDA on tensor data, several methodologies
have been introduced in the past. Among them, we can cite discriminant
analysis with tensor representation (DATER) \cite{Yan1}, Tensor
subspace analysis (TSA) \cite{He}, multilinear discriminant analysis
(MDA)  \cite{Yan,Yan1} and constrained multilinear discriminant analysis
(CMDA)  \cite{li}. The principle of MDA is to find a lower dimensional
tensor subspace represented by  orthonormal matrices. However, the main drawback of the current approaches relies on its optimization step which is based on heuristic optimization
approaches.

In this paper, we propose a new way to compute linear discriminant
analysis from three-order tensors. This new method is based on recent
developments on tensor-tensor products \cite{Kilmer2011,Kernfeld2015}.
The first work on this issue is due to Kilmer et al. \cite{Kilmer2011}  where they introduced 
the notion of the t-product which allows to mutiply easily 3-order
tensors. The multiplication uses a convolution-type operation which can be advantageously computed by Fast Fourier Transform (FFT). More recently, Kernfeld et
al. extended this approach and defended the principle that any tensor-tensor
product can be defined with arbitrary invertible linear transforms \cite{Kernfeld2015}. As an example, they introduced the tensor cosine transform product which is an  alternative of the t-product and  can be efficiently computed by the Discret Cosine Transform (DCT). 

Motivated by these works, we propose to solve Tensor Linear Discriminant
Analysis (TLDA) by using the concept of "transform-domain" and the cited new 
family of tensor-tensor products. We show here that moving TLDA to
the transform domain (TD) makes its resolution easier and more natural.
First, one of the fundamental step of the previous tensor LDA approaches
is to compute some tensor-matrix products that  involve a succession
of tensor unfoldings. In this context, the corresponding optimization
problem relies on the determination of a set of projection matrices i.e two-dimensional projective subspaces. Such an approach shows the
difficulty to formulate this problem in fully tensor form by directly
searching for a projective tensor instead of a set of projective matrices.
By moving in the transform domain, the main ingredients of TLDA can
be formulated as tensor-tensor products which are efficiently computed
by a sequence of matrix-matrix products. Secondly, a key point of this approach is that the optimization problem can be solved in the transform domain benefiting from its properties.  More precisely, each frontal face of the transformed tensor are processed independently to build independent optimization sub-problems easier to address. The obtained solutions are returned in the form of a projective  tensor by inverse transform.  
This paper is organized as follows: Section 2 introduces the notation and the main definitions. Section 3 recalls the main definitions of the tensor-tensor product
under the concept of invertible linear transforms. In Section 4, after recalling the principle of the standard LDA (matricial case) and the tensor LDA (TLDA) formulated with the n-product, we present a new multilinear discriminant analysis based on third-order tensors and formulated with new tensor-tensor products. Section 5 analyses and compares the performance of the proposed approach with recent tensor based LDA approaches on several multidimensional data sets. Section 6  ends the paper with a conclusion.

\section{Notation and preliminaries}

Scalars, vectors, matrices and high-order tensors will be denoted by
lowercase letters, e.g. $a$, boldface lowercase letters, e.g. $\boldsymbol{a}$,
capital letters, e.g. $A$ and Euler script letters, e.g. $\mathcal{A}$,
respectively. In this work, we will limit our study to third-order
tensors. Third-order tensors are compact and well adapted to represent
multidimensional data from vision based applications such as face
identification, video monitoring or classification of multispectral
images. Let $\mathcal{A}\in\mathbb{R}^{n_{1}\times n_{2}\times n_{3}}$
be a third-order tensor. By convention, the first dimension is devoted
to the \textit{pixels} of the images, the second dimension for the
\textit{number of images} and the third dimension is for the \textit{number
	of modalities of the images}. As an illustrative example, if we consider
a sequence of $l$ color images of size $n\times m$, the corresponding
third order tensor will be sized as follows: $n_{1}=nm$, $n_{2}=l$
and $n_{3}=3$. \\
The tensor $\mathcal{A}$ is sampled by a triplet of indexes $(i,j,k)$ which
allows to select different subparts of $\mathcal{A}$. By fixing the
whole set of indexes we obtain a scalar entry of $\mathcal{A}$ denoted
by $a_{ijk}$ and  fixing two indexes over three, we select a \textit{
	fiber} of $\mathcal{A}$. We will denote a column (\textit{1-mode}),
row (\textit{2-mode}) and tube (\textit{3-mode}) fiber by $\boldsymbol{a}_{.jk}$,
$\boldsymbol{a}_{i.k}$ and $\boldsymbol{a}_{ij.}$, respectively.
Lastly, by fixing one index over three, we define a \textit{slice}
of $\mathcal{A}$. Therefore, slices  are declined in three modes: horizontal (\textit{1-mode}),
lateral (\textit{2-mode}) and frontal (\textit{3-mode}) slides which
are represented by $A_{i::}$, $A_{:j:}$ and $A_{::k}$, respectively.
In the sequel, the $k^{th}$ frontal slide of a third-order tensor
will be denoted more compactly by$A^{(k)}$. \\
Manipulating tensors needs specific algebra. Here we just list some
definitions which are directly relevant to this paper. For a detailed
description, see  for example \cite{Kilmer2011,kolda1}.
 .
\subsection{Tensor unfolding}
Tensor unfolding or flattening consists in reordering the elements
of a tensor into a matrix. Consider the general case of a $N^{th}$order
tensor $\mathcal{A}\in\mathbb{R}^{n_{1}\times n_{2}\times...\times n_{N}}$,
flattening $\mathcal{A}$ along the $k^{th}$ mode or the \textbf{\textit{k-mode}}
matricization of $\mathcal{A}$ gives a matrix denoted $A_{(k)}$
which consists in arranging the \textit{k-mode} fibers to be the columns
of the resulting matrix. 

\subsection{Tensor products}
Let us recall  several tensor products.
\begin{definition}[\textbf{k-mode product}]\label{defkmode}
	Consider the general case of na $N^{th}$ order tensor $\mathcal{A}\in\mathbb{R}^{n_{1}\times n_{2}\times...\times n_{N}}$.
	The \textbf{\textit{k-mode}} product of $\mathcal{A}$ with a matrix
	$U\in\mathbb{R}^{m\times n_{k}}$ is a new tensor $\mathcal{B}\in\mathbb{R}^{n_{1}\times...\times n_{k-1}\times\boldsymbol{m}\times n_{k+1}\times...\times n_{N}}$
	defined by 
	\begin{equation}
		\mathcal{B}=\mathcal{A}\times_{k}U\label{eq:A_1}
	\end{equation}
	which is equivalent to the following matrix-matrix product 
	\begin{equation}
		B_{(k)}=UA_{(k)}\label{eq:A_2}
	\end{equation}
		where $A_{(k)}$ and $B_{(k)}$ denotes the \textit{k-mode} matricization
	(see section 2.1) of $\mathcal{A}$ and $\mathcal{B}$, respectively. 
\end{definition}
\\

\begin{definition}[\textbf{generalization of the k-mode product}]
	Consider the general case of an $N^{th}$order tensor $\mathcal{A}\in\mathbb{R}^{n_{1}\times n_{2}\times...\times n_{N}}$.
	The multiplication of $\mathcal{A}$ with a set of $N$ matrices $\left\{ U_{k}\in\mathbb{R}^{m_{k}\times n_{k}},k=1,...,N\right\} $
	is defined by
\begin{equation}
	\mathcal{B}=\mathcal{A}\Pi_{i=1}^{n}\times_{i}U_{i}=\mathcal{A}\times_{1}U_{1}\times_{2}U_{2}\times...\times_{N}U_{N}. \label{eq:A_3}
\end{equation} 
The \textit{k-mode} matricization of $\mathcal{B}$ can be obtained by
\begin{equation}
	B_{(k)}=U_{k}A_{(k)}U_{\bar{k}}^{\top}\label{eq:A_4}
\end{equation}
where $U_{\bar{k}}=U_{1}\otimes...\otimes U_{k-1}\otimes U_{k+1}\otimes...\otimes U_{N}$.
\end{definition}
\\
\begin{definition}[\textbf{face-wise product} \cite{Kernfeld2015}]
Let $\mathcal{A}\in\mathbb{R}^{n_{1}\times n_{2}\times n_{3}}$ and $\mathcal{B}\in\mathbb{R}^{n_{2}\times m\times n_{3}}$ be two third-order tensors, then the \textbf{\textit{face-wise}} product between them consists in computing a matrix-matrix product between the 3-mode slides
	of $\mathcal{A}$ and $\mathcal{B}$ as follows 
	\begin{equation}
		\left(A\triangle B\right)^{(i)}=A^{(i)}B^{(i)},\quad i=1,...,n_{3}\label{eq:B}
	\end{equation}
\end{definition}

\subsection{Specific block matrices}
Tensor-tensor products require the definition of specific structured
block matrices build from the frontal slices of the third-order tensor.
We recall here some definitions 

\begin{definition}
The \textbf{\textit{Toeplitz-plus-Hankel}} matrix of the tensor $\mathcal{A}$
is a $n_{1}n_{3}\times n_{2}n_{3}$ block matrix composed of the frontal
slices $A^{(i)}$, $i=1,\ldots,n_3$ of $\mathcal{A}$ and defined by
\begin{equation}
	\textrm{mat}(\mathcal{A})=\left(\begin{array}{cccc}
		A^{(1)} & A^{(2)} & \cdots & A^{(n_{3})}\\
		A^{(2)} & A^{(1)} & \cdots & A^{(n_{3}-1)}\\
		\vdots & \vdots & \cdots & \vdots\\
		A^{(n_{3})} & A^{(n_{3}-1)} & \cdots & A^{(1)}
	\end{array}\right)+\left(\begin{array}{cccc}
		A^{(1)} & \cdots & A^{(n_{3})} & \boldsymbol{0}\\
		\vdots & \iddots & \iddots & A^{(n_{3})}\\
		A^{(n_{3})} & \boldsymbol{0} & \iddots & \vdots\\
		\boldsymbol{0} & A^{(n_{3})} & \cdots & A^{(1)}
	\end{array}\right)\label{eq:C}
\end{equation}
where $\boldsymbol{0}$ denotes the zero matrix of size $n_{1}\times n_{2}$. We
will denote $ten$ the inverse operator such as: $$\textrm{ten(mat}(\mathcal{A}))=\mathcal{A}.$$
\end{definition}

\begin{definition}
The \textbf{block circulant} matrix of a tensor $\mathcal{A}$ is
the  $n_{1}n_{3}\times n_{2}n_{3}$ block matrix composed by  the frontal
slices of $\mathcal{A}$ and defined by 
\begin{equation}
	\textrm{bcirc}(\mathcal{A})=\left(\begin{array}{cccc}
		A^{(1)} & A^{(n_{3})} & \cdots & A^{(2)}\\
		A^{(2)} & A^{(1)} & \cdots & A^{(3)}\\
		\vdots & \vdots & \cdots & \vdots\\
		A^{(n_{3})} & A^{(n_{3}-1)} & \cdots & A^{(1)}
	\end{array}\right)\label{eq:D}
\end{equation}
\end{definition}

\begin{definition}
The \textbf{block diagonal} matrix of a tensor $\mathcal{A}$ is the  $n_{1}n_{3}\times n_{2}n_{3}$
block matrix composed by  the frontal slices of $\mathcal{A}$ and
defined by  
\begin{equation}
	\textrm{bdiag}\left(\mathcal{A}\right)=\left(\begin{array}{cccc}
		A^{(1)} & 0 & \cdots & 0\\
		0 & A^{(2)} & \cdots & 0\\
		\vdots & \vdots & \ddots & \vdots\\
		0 & 0 & \cdots & A^{(n_{3})}
	\end{array}\right)\label{eq:E}
\end{equation}
\end{definition}

\section{Tensor-tensor products with invertible linear transform}
Recently, a new type of tensor-tensor products, called t-product, has
been proposed in \cite{kilmer0}. The t-product generalizes matrix multiplication
for third-order tensors. It is based on a convolution-like operation
which is efficiently computed by the Fast Fourier Transform (FFT).
This work opens the way towards the idea that there exits a \textit{transform
domain} where the tensor-tensor product can be defined. Motivated
by this idea, Kernfeld et al. \cite{Kernfeld2015}   extend this concept by introducing
a new family of tensor-tensor products which can be efficiently computed
in a \textit{transform domain } for any invertible linear transform.
To illustrate their principle, they defined the c-product which is
an alternative to the t-product which can be efficiently computed
in the transform domain via the discret cosinus transform (DCT). In the
sequel, we will consider that the result of the transformation is
at more of complex type.
We first  recall the main properties and definitions introduced
in  \cite{Kernfeld2015} .\\
\begin{definition}
Let $L:\mathbb{R}^{1\times1\times n_{3}}\rightarrow\mathbb{C}^{1\times1\times n_{3}}$
be an invertible transform and $\mathcal{A}\in\mathbb{R}^{n_{1}\times n_{2}\times n_{3}}$
be a third-order tensor. $L$ transforms any tube fibers $\boldsymbol{a}\in\mathbb{R}^{1\times1\times n_{3}}$
of $\mathcal{A}$ into $\boldsymbol{\widetilde{a}}\in\mathbb{\mathbb{C}}^{1\times1\times n_{3}}$
in the following way
\[
\boldsymbol{\widetilde{a}}^{(k)}=\left(L(\boldsymbol{a})\right)^{(k)}=\left(M.\textrm{vec}(\boldsymbol{a})\right)_{k}\quad k=1...n_{3}
\]
where $M$ is a $n_{3}\times n_{3}$ invertible matrix associated
to $L$ and $\textrm{vec}(\boldsymbol{a})$ is the vector in  $\mathbb{R}^{n_3}$ whose elements are the elements of the tube $\boldsymbol{a}$ . From a practical point of view, $\mathcal{\widetilde{A}}\in\mathbb{\mathbb{C}}^{n_{1}\times n_{2}\times n_{3}}$,
the transform domain version of $\mathcal{A}$, can be efficiently
computed as follows
\begin{equation}
	\mathcal{\widetilde{A}}=L(\mathcal{A})=\mathcal{A}\times_{3}M. \label{eq:B1}
\end{equation}
Similarly, we have 
\begin{equation}
	\mathcal{\mathcal{A}}=L^{-1}(\widetilde{\mathcal{A}})=\widetilde{\mathcal{A}}\times_{3}M^{-1}\label{eq:B2}
\end{equation}
where $\times_{3}$ is the \textbf{\textit{3-mode}} product as defined
in (\ref{defkmode}).
\end{definition}

\noindent Notice  that the matrix M is specific to the transform $L$
and to the corresponding tensor-tensor product (see appendix 6.1).  \\

\begin{definition}
Let $\textrm{op}_{L}(\mathcal{A})$ be a structured block matrix build
from $\mathcal{A}$ and specific to the L-transform (see appendix 6.2). Let $\boldsymbol{A}\in\mathbb{\mathbb{C}}^{n_{1}n_{3}\times n_{2}n_{3}}$
be a block diagonal matrix where each block represents a frontal slice
of the transform tensor $\mathcal{\widetilde{A}}$, then it can be
shown that $\widetilde{A}$ results from the block diagonalization
of $\textrm{op}_{L}(\mathcal{A}$) as
\begin{equation}
	\boldsymbol{A}=bdiag\left(\widetilde{\mathcal{A}}\right)=(M\otimes I_{n_{1}})\textrm{op}_{L}(\mathcal{A})(M^{-1}\otimes I_{n_{2}})\label{eq:B7-1}
\end{equation}
where $M$ is the $n_{3}\times n_{3}$ transform matrix associated
to $L$ and $I_{n}$ iis the $n\times n$ identity matrix.
\end{definition}

\medskip
\noindent The $L$-product of   two tensors is defined as follows
\begin{definition}
Let $*_{L}:\mathbb{R}^{m\times l\times n_{3}}\times\mathbb{R}^{l\times p\times n_{3}}\rightarrow\mathbb{R}^{m\times p\times n_{3}}$
be the product operator in $L$ defined such as 
\begin{equation}
	L\left(\mathcal{A}*_{L}\mathcal{B}\right)=L\left(\mathcal{A}\right)\bigtriangleup L\left(\mathcal{B}\right)\label{eq:B2-1}
\end{equation}
where $\bigtriangleup$ is the \textbf{\textit{face-wise}} product
as defined in (\ref{eq:B}). Let $\mathcal{C}\in\mathbb{R}^{m\times p\times n_{3}}$
be the result of the $L$-product between $\mathcal{A}$ and $\mathcal{B}$
, then we have 
\begin{equation}
	\mathcal{C}=\mathcal{A}*_{L}\mathcal{B}=L^{-1}\left(L\left(\mathcal{A}\right)\bigtriangleup L\left(\mathcal{B}\right)\right)\label{eq:B2-2}
\end{equation}
\end{definition}

\noindent The main known and used  $L$-products are the t-product and the c-product. For these two products, the matrix $M$ is given as follows. \\
For the c-product whicgh is based on DCT, $M$ is given by 
\begin{equation}
	M=W^{-1}C\left(I+Z)\right)\label{eq:B3}
\end{equation}
where $C$ is the $n_{3}\times n_{3}$ DCT matrix where each entry
is defined by
\begin{equation}
	c_{ij}=\sqrt{\frac{2-\delta_{ij}}{n_{3}}}cos\left(\frac{(i-1)(2j-1)}{2n_{3}}\right)\;\left(i,j\right)=(1,...,n_{3})\label{eq:B4}
\end{equation}
where $\delta$ is the Kronecker indicator. $W=diag(\boldsymbol{c}_{.1})$
is the diagonal matrix build from the first column of $C$ and $Z$
is an $n_{3}\times n_{3}$ circulant upshift matrix.

\noindent For the t-product which is based on FFT, $M$ is given by 
\begin{equation}
	M=F\label{eq:B5}
\end{equation}
where $F$ is the $n_{3}\times n_{3}$ FFT matrix where each entry
is defined by
\begin{equation}
	f_{ij}=exp\left(-j2\pi\frac{(i-1)(j-1)}{n_{3}}\right)\label{eq:B6}
\end{equation}
We also notice that  $\textrm{op}(\mathcal{A})$
is defined by
\[
\textrm{op}_L(\mathcal{A})=\left\{ \begin{array}{ccc}
	\textrm{bcirc}(\mathcal{A}) &  & \textrm{for the t-product}\\
	\textrm{mat}(A) &  & \textrm{for the c-product}
\end{array}\right.
\]
where the operators \textit{bcirc} and \textit{mat }are defined in
Section 1.

\noindent From the relation (\ref{eq:B7-1}), it can be shown that
$\mathcal{C}=\mathcal{A}*_{L}\mathcal{B}$ is equivalent to compute
$\boldsymbol{C}=\boldsymbol{A}\boldsymbol{B}$ in the transform domain. 
Algorithm \ref{lprod} allows us to compute in an efficient way the $L$-product of the tensors  $\mathcal{A}$ and $\mathcal{B}$.
\newpage
\begin{algorithm}[!h] 	
	\caption{Tensor–tensor product via the operator $L$}\label{lprod}
	\begin{algorithmic} 
		\STATE \textbf{Inputs}: $\mathcal{A}\in\mathbb{R}^{m\times l\times n_{3}}$, $\mathcal{B}\in\mathbb{R}^{l\times p\times n_{3}}$ \\
		\STATE \textbf{Output}: $\mathcal{C}\in\mathbb{R}^{m\times p\times n_{3}}$\\
		\STATE $\widetilde{\mathcal{A}}=L(\mathcal{A})$
	\STATE $\widetilde{\mathcal{B}}=L(\mathcal{B})$
		\FOR{$i=1,\ldots,n_{3}$} 
		\STATE $\mathcal{C} ^{(i)}=\mathcal{A} ^{(i)}\mathcal{B} ^{(i)}$.
		\ENDFOR \\
		$\mathcal{C}=L^{-1}(\widetilde{\mathcal{C}})$
	\end{algorithmic}
\end{algorithm}
Some basic algebraic properties are associated to ${L}$-product  such as
associativity, distribution over addition and invertibility.
$*_{L}$ has a identy element, the hermitian
transpose, norms and inner products (for more details see \cite{Kernfeld2015}). 

\section{Tensor LDA using the L-product: $*_{L}$- TLDA }
Linear discriminant analysis is a supervised dimensionality reduction
method which aims to find a low-dimensional projective subspace which
best separates $n$ training data vectors $\boldsymbol{x}_{1},\boldsymbol{x}_{2},...,\boldsymbol{x}_{n}$
into $c$ classes or clusters. In the sequel, we will consider that
each data vector $\boldsymbol{x}_{k}$ belongs to a class indexed
as $l_{k}\in\left\{ 1,2,...,c\right\} $ and each class $i$ is defined
by a set of indices $C_{i}$ of length $n_{i}$ such as $n=\sum_{i=1}^{c}n_{i}$
.
In this section, considering the previous notation, we propose to develop the tensor linear discriminant analysis using the $L$-product  which will be denoted by  $*_{L}$-TLDA. Before
all, let us recall the formulation of LDA in the matricial case and
the tensor LDA (using the $n$-product).

\subsection{LDA}
Consider that the training samples are collected into a matrix $X=\left[\boldsymbol{x}_{1},\boldsymbol{x}_{2},...,\boldsymbol{x}_{n}\right]$ $\in\mathbb{R}^{n_{1}\times n}$
where each component $\boldsymbol{x}_{i}$ is an $n_{1}$ dimensional
data vector. Let $\boldsymbol{m}=\frac{1}{n}\sum_{i=1}^{n}\boldsymbol{x}_{i}$
be the global centroid of $X$ and $\boldsymbol{m}_{i}=\frac{1}{n_{i}}\sum_{i\in C_{i}}\boldsymbol{x}_{i}$
be the centroid of the data vectors belonging to the cluster $i$,
then the goal of LDA can be defined as follows \\
\begin{definition}[\textbf{LDA}]
	Let $V=\left[\boldsymbol{v}_{1},\boldsymbol{v}_{2},...,\boldsymbol{v}_{m}\right]\in\mathbb{R}^{n_{1}\times m}$
	be a matrix defining a low-dimensional projective subspace ($m\ll n_{1}$)
	and 
	\begin{equation}
		\left\{ \begin{array}{clc}
			\psi_{B}(V)= & \sum_{j=1}^{c}n_{j}\left\Vert V^{T}\left(\boldsymbol{m}_{j}-\boldsymbol{m}\right)\right\Vert _{F}^{2} & \quad(a)\\
			\psi_{W}(V)= & \sum_{j=1}^{c}\sum_{i\in C_{j}}\left\Vert V^{T}\left(\boldsymbol{x}_{i}-\boldsymbol{m}_{j}\right)\right\Vert _{F}^{2} & \quad(b)
		\end{array}\right.\label{eq:LDA_1}
	\end{equation}
	be the between and within scatters measured in the projective subspace
	$V$, respectively. Then, LDA consists in finding the projective subspace
	$V^{*}$ which maximizes the between scatter measure (a) while minimizing
	the within scatter measure (b) which can be formulated as follows
		\begin{equation}
		V^{*}=\underset{V^{\top}V=I}{max}\quad\frac{\psi_{B}(V)}{\psi_{W}(V)}\label{eq:LDA_2}
	\end{equation}
\end{definition}
Let us introduce 
\begin{equation}
	\left\{ \begin{array}{clc}
		S_{B}= & \sum_{j=1}^{c}n_{j}\left(\boldsymbol{m}_{j}-\boldsymbol{m}\right)\left(\boldsymbol{m}_{j}-\boldsymbol{m}\right)^{T} & \quad(a)\\
		S_{W}= & \sum_{j=1}^{c}\sum_{i\in C_{j}}\left(\boldsymbol{x}_{i}-\boldsymbol{m}_{j}\right)\left(\boldsymbol{x}_{i}-\boldsymbol{m}_{j}\right)^{\top} & \quad(b)
	\end{array}\right.\label{eq:LDA_11}
\end{equation}
the between and within scatter matrices, respectively, then it can be shown easily that the problem (\ref{eq:LDA_2}) can be re-written as 
\begin{equation}
	V^{*}=\underset{V^{\top}V=I}{max}\quad\frac{\mathrm{Tr}\left(V^{\top}S_{B}V\right)}{\mathrm{Tr}(V^{\top}S_{W}V)}\label{eq:LDA_3}
\end{equation}
Problem (\ref{eq:LDA_3}), also referred as the \textbf{trace ratio}
problem, is non convex and does not have a closed-form solution. Fortunately,
it can be shown that it is equivalent to a trace difference problem
\begin{equation}
	V^{*}=\underset{V^{\top}V=I}{max}\quad\mathrm{Tr} (V^T    (S_B   -\rho    S_W) V)\label{eq:LDA_31}.
\end{equation}
which can be solved iteratively by the Newton-Lanczos
algorithm {\cite{ngobellalij}}. Algorithm \ref{Saadalgo} summarizes the main steps of the maximization of the trace ratio problem with the Newton-Lanczos algorithm (\ref{eq:LDA_31}).
\begin{algorithm}[h!]
	\caption{Newton-Lanczos algorithm for Trace Ratio}	\label{Saadalgo}
	
	{\bf Input} :  two  matrix $S_A$ and $S_B$.
	\begin{itemize}
		\item Select a unitary matrix $V$ with $k$ columns and compute 
		$\rho=\displaystyle \frac{{\text Trace}(V^TS_BV)}{{\text Trace}(V^TS_WV)}$.
		\item Until convergence do:
		\begin{enumerate}
			\item Call the Lanczos algorithm to compute the largest $k$ eigenvalues $\lambda_1(\rho),\ldots,\lambda_k (\rho),$  of $S_B-\rho S_W$ and the associated eigenvectors: $V=[v_1,\ldots,v_k]$.
			\item Set $\rho=\displaystyle \frac{{\text Trace}(V_k^TS_BV_k)}{{\text Trace}(V_k^TS_WV_k)}$ and go to Step 1.
		\end{enumerate}
		\item	EndDo
	\end{itemize}	
\end{algorithm}

\noindent It has been shown in \cite{ngobellalij} that this algorithm converges to a global optimum. However, the drawback of this procedure is the repeated calls to an eigensolver which can be time-consuming when the dimensionality of the data is very large. Another critical point concerns the choice of the reduced dimension. Indeed, the output dimension $m$ is bounded by the rank of the matrix $S_W$ since the $rank(S_B   -\rho    S_W)=rank(S_W)$;  ($rank(S_B)<rank(S_W)$). Since $rank(S_W)<n-c$, then $m$ is at most $n-c$. As a consequence, the optimal output dimension is related to the sample size $n$ and its selection may be also time-consuming when the size of the training sample is very large. 
These observations often lead to replace the \textbf{trace ratio} problem by the
simpler, but not equivalent \textbf{ratio trace} problem
\begin{equation}
	V^{*}=\underset{V^{\top}V=I}{max}\quad\mathrm{Tr}\left((V^{\top}S_{W}V)^{-1}(V^{\top}S_{B}V)\right)\label{eq:LDA_4}
\end{equation}
which has a closed-form solution. It is equivalent to solve the following generalized eigenvalue problem
\begin{equation}
	S_{B}U=\Lambda S_{W}U\label{eq:LDA_5}
\end{equation}
where $U$ denotes the matrix of eigenvectors and $\Lambda$ the diagonal
matrix of eigenvalues. Thus the projection matrix $V$ is explicitly
characterized through the eigen-decomposition of the matrix $S_{W}^{-1}S_{B}$
if $S_{W}$ is nonsingular. Moreover, the dimension of the projective
subspace is defined by the rank of $S_{B}$ which implies that $m$
is at most $c-1$;  ($m=c-1$ when data are linearly independent). When
$S_{W}$ becomes singular, the problem is said ``undersampled'',
i.e the sample size is smaller than the dimension of the data. A common strategy is to introduce regularization into the problem
(\ref{eq:LDA_5}) which translates into
\begin{equation}
	(S_{W}+\gamma I)^{-1}S_{B}U=\Lambda U\label{eq:LDA_6}
\end{equation}
where $I$ is the identity matrix and $\gamma>0$ the regularization
parameter. The value of $\gamma$ must be chosen with care and its
selection can be  obtained by cross validation. 

\subsection{The $\times_{n}$-TLDA}
Consider the general case where each data sample is represented by
an $N^{th}$ order tensor $\mathcal{X}_{i}\in\mathbb{R}^{n_{1}\times n_{2}\times...\times n_{N}}$
and the sample set by a $(N+1)^{th}$ order tensor $\mathcal{X}\in\mathbb{R}^{n_{1}\times n_{2}\times...\times n_{N}\times\boldsymbol{n}}$.
Let $\mathcal{M}=\frac{1}{N}\sum_{i=1}^{n}\mathcal{X}_{i}$ be the
global mean of $\mathcal{X}$ and $\mathcal{M}_{i}=\frac{1}{n_{i}}\sum_{i\in C_{i}}\mathcal{X}_{i}$
be the mean of the cluster $i$. Then the goal of the tensor LDA or
multilinear DA can be defined as follows
\medskip
\begin{definition}[$\times_{n}$-\textbf{TLDA}]
	Let $V_{k}\mid_{k=1}^{K}$ be a set low-dimensional projective matrices
	of size $n_{k}\times m_{k}$ with $m_{k}\ll n_{k}$ and let
	\begin{equation}
		\left\{ \begin{array}{ccc}
			\psi_{B}(V_{k}\mid_{k=1}^{K})= & \sum_{j=1}^{c}n_{j}\left\Vert \left(\mathcal{M}_{j}-\mathcal{M}\right)\Pi_{k=1}^{K}\times_{k}V_{k}\right\Vert _{F}^{2} & \quad(a)\\
			\psi_{W}(V_{k}\mid_{k=1}^{K})= & \sum_{j=1}^{c}\sum_{i\in C_{j}}\left\Vert \left(\mathcal{X}_{i}-\mathcal{M}_{j}\right)\Pi_{k=1}^{K}\times_{k}V_{k}\right\Vert _{F}^{2} & \quad(b)
		\end{array}\right.\label{eq:nTDA_1}
	\end{equation}
		be the between and within scatters measured in the set of projective
	subspace $U_{k}\mid_{k=1}^{K}$, respectively. Then, $\times_{n}$-TLDA
	consists in finding a set projective subspace $V^{*}$ that maximizes
	the between scatter measure (a) while minimizing the within scatter
	measure (b), i.e. 
	\begin{equation}
		V_{k}^{*}=\underset{V_{k}}{max}\quad\frac{\psi_{B}(V_{k})}{\psi_{W}(V_{k})},\;\; k=1,\ldots,N.\label{eq:nTDA_2}
	\end{equation}
\end{definition}
However, the objective function (\ref{eq:nTDA_2}) has no closed-form
solution due to that the $V_{k}$s, $k=1,\ldots,N$  depends on each other
and the standard procedure is to solve it by an iterative optimization
procedure. Considering that $\left\Vert \mathcal{X}\right\Vert =\left\Vert X_{(k)}\right\Vert _{F}$
and $\left\Vert X\right\Vert ^{2}=\mathrm{Tr}\left(X^{\top}X\right)=\mathrm{Tr}\left(XX^{\top}\right)$,
if we assume that $K-1$ projective matrices $V_{i}\mid_{i=1,i\neq k}^{K}$
have been previously computed, then $V_{k}$ is updated by maximizing 
\begin{equation}
	V_{k}^{*}=\underset{V_{k}}{max}\quad\frac{\mathrm{Tr}\left(V_{k}^{\top}S_{B(k)}V_{k}\right)}{\mathrm{Tr}\left(V_{k}^{\top}S_{W(k)}V_{k}\right)}\label{eq:nTDA_3}
\end{equation}
where $S_{B(k)}$ and $S_{W(k)}$ denote the between-class and within-class
scatter matrices along the $k^{th}$ mode, respectively and defined by
\begin{equation}
	\left\{ \begin{array}{clc}
		S_{B(k)}= & {\displaystyle \sum_{j=1}^{c}}n_{j}\left((\mathcal{M}_{j}-\mathcal{M})_{(k)}\right)V_{\bar{k}}^{\top}V_{\bar{k}}\left((\mathcal{M}_{j}-\mathcal{M})_{(k)}\right)^{\top} & \quad(a)\\
		S_{W(k)}= & \sum_{j=1}^{c}\sum_{i\in C_{j}}\left((\mathcal{X}_{i}-\mathcal{M}_{j})_{(k)}\right)V_{\bar{k}}^{\top}V_{\bar{k}}\left((\mathcal{X}_{i}-\mathcal{M}_{j})_{(k)}\right)^{\top} & \quad(b)
	\end{array}\right.\label{eq:nTDA_4}
\end{equation}
where $V_{\bar{k}}=V_{K}\otimes...\otimes V_{k+1}\otimes V_{k}\otimes...\otimes V_{1}$
(see definition 2, (\ref{eq:A_4})) and the terms $(\mathcal{A}-\mathcal{B})_{(k)}$
in (\ref{eq:nTDA_4}) denote the $k$-mode matricization of the tensor
$\mathcal{A}-\mathcal{B}$ (see section 2.1). This iterative optimization
procedure, also called $k$-mode optimization, have been originaly
introduced in \cite{Yan} and became the central part of several
work to solve multilinear discriminant analysis (MDA), \cite{li,Yan,Yan1,Lu2009}.
However, all these methods solve the MDA problem from heuristic optimization
procedures that do not rigorously optimize the MDA objective. In the
sequel, we propose to solve MDA objective using the  $L$-tensor-tensor
products. The corresponding optimization problem can be moved into
an invertible transform domain in which a closed-form solution exits.
In the sequel, we propose to develop TLDA using the $L$-product,i.e. $*_{L}$-TLDA.

\subsection{The $*_{L}$-TLDA}
Assume the learning data set is composed of $n$ samples and each sample
is represented by a third-order tensor, i.e. $\left\{ \mathcal{X}_{i}\in\mathbb{R}^{n_{1}\times1\times n_{3}},i=1,...,n\right\} $. The sample set can by represented by a unique third-order tensor
$\mathcal{X}\in\mathbb{R}^{n_{1}\times\boldsymbol{n}\times n_{3}}$.
Let $\mathcal{M}=\frac{1}{n}\sum_{i=1}^{n}\mathcal{X}_{i}$ be the
global centroid of $\mathcal{X}$ and $\mathcal{M}_{i}=\frac{1}{n_{i}}\sum_{i\in C_{i}}\mathcal{X}_{i}$
be the centroid of the tensors belonging to the cluster $i$ , then the goal of $*_{L}$-TLDA can be defined as follows\\

\begin{definition}[$*_{L}$-\textbf{TLDA}]
Let $\mathcal{V}\in\mathbb{R}^{n_{1}\times K\times n_{3}}=\left[\mathcal{V}_{1},\mathcal{V}_{2},...,\mathcal{V}_{K}\right]$
, $\mathcal{V}_{i}\in\mathbb{R}^{n_{1}\times1\times n_{3}}$ be a
projective third-order tensor where $K$ denotes the dimension of
the projective subspace and 
\begin{equation}
	\left\{ \begin{array}{cl}
		\psi_{B}(\mathcal{V})= & {\displaystyle \sum_{j=1}^{c}}\left\Vert \mathcal{V}^{T}*_{L}\left(\mathcal{M}_{j}-\mathcal{M}\right)\right\Vert _{F}^{2}\\
		\psi_{W}(\mathcal{V})= & \sum_{j=1}^{c}\sum_{i\in N_{j}}\left\Vert \mathcal{V}^{T}*_{L}\left(\mathcal{X}_{i}-\mathcal{M}_{j}\right)\right\Vert _{F}^{2}
	\end{array}\right. \label{eq:scatter}
\end{equation}
where $\psi_{B}(\mathcal{V})$ and $\psi_{W}(\mathcal{V})$ denote
the \textit{between} and \textit{within} scatter measures, respectively.
Then, the goal of $*_{L}$- LDA is to find a projective third-order
tensor $\mathcal{V}\in\mathbb{R}^{n_{1}\times k\times n_{3}}$ which
maximizes the following objective function
\begin{equation}
	\mathcal{V}^{*}\in\mathbb{R}^{n_{1}\times K\times n_{3}}:\underset{\mathcal{V}}{max}\;\frac{\psi_{B}(\mathcal{V})}{\psi_{W}(\mathcal{V})}\label{eq:C-1}
\end{equation}
\end{definition}
Notice  that the product $\mathcal{V}^{T}*_{L}\left(\mathcal{M}_{i}-\mathcal{M}\right)$
represents the orthogonal projection of one lateral slice in $\mathbb{R}^{n_{1}\times1\times n_{3}}$ or
$n_{1}$ tube fibers onto $\mathcal{V}$ which generate $K$ tube
fibers. The same remark can be made for the product $\mathcal{V}^{T}*_{L}\left(\mathcal{X}_{i}-\mathcal{M}_{i}\right)$.
Second, (\ref{eq:C-1}) can be solved more easily in the transform
domain by using the fact that the $L$-product between tensors can
be advantageously replaced by simple matrix-matrix products between
the transform versions of the tensors (\ref{eq:B2-1}). Then, the
$L$-product is recovered by inverse transform. Let us develop the
main steps:\\
From Definition \ref{eq:B7-1}, we denote $\boldsymbol{V}=bdiag(\widetilde{\mathcal{V}})$,
$\boldsymbol{M}=bdiag(\widetilde{\mathcal{M}})$, $\boldsymbol{M}_{i}=bdiag(\mathcal{\widetilde{M}}_{i})$
and $\boldsymbol{X}_{i}=bdiag(\mathcal{\widetilde{X}}_{i})$ the block
diagonal matrices build from the transform versions of the tensors $\mathcal{V}$, $\mathcal{M}$, $\mathcal{M}_{i}$ and $\mathcal{X}_{i}$, respectively. Since $\mathcal{C}=\mathcal{A}*_{L}\mathcal{B}$ is
equivalent to compute $\boldsymbol{C}=\boldsymbol{A}\boldsymbol{B}$
in the transform domain, the $L$-products $\mathcal{V}^{T}*_{L}\left(\mathcal{M}_{i}-\mathcal{M}\right)$
and $\mathcal{V}^{T}*_{L}\left(\mathcal{X}_{i}-\mathcal{M}_{i}\right)$
can be computed via a simple product of block diagonal matrices, i.e.
$\boldsymbol{V}^{\top}\left(\boldsymbol{M}_{i}-\boldsymbol{M}\right)$
and $\boldsymbol{V}^{\top}\left(\boldsymbol{X}_{i}-\boldsymbol{M_{i}}\right)$.
Then, the computation of $\psi_{B}(\mathcal{V})$ becomes in the transform
domain
\begin{equation}
	\begin{array}{rl}
		\widetilde{\psi}_{B}(\boldsymbol{V})= & \sum_{j=1}^{c}\left\langle \boldsymbol{V}^{\top}\left(\boldsymbol{M}_{j}-\boldsymbol{M}\right),\boldsymbol{V}^{\top}\left(\boldsymbol{M}_{j}-\boldsymbol{M}\right)\right\rangle \\
		= & \sum_{j=1}^{C}\textrm{Trace}\left(\boldsymbol{V}^{\top}\left(\boldsymbol{M}_{j}-\boldsymbol{M}\right)\left(\boldsymbol{M}_{j}-\boldsymbol{M}\right)^{\top}\boldsymbol{V}\right)\\
		= & \mathrm{\mathrm{Trace}}\left(\boldsymbol{V}^{\top}\sum_{j=1}^{c}\left(\boldsymbol{M}_{j}-\boldsymbol{M}\right)\left(\boldsymbol{M}_{j}-\boldsymbol{M}\right)^{\top}\boldsymbol{V}\right)\\
		= & \mathrm{Trace}\left(\boldsymbol{V}^{\top}\boldsymbol{S}_{B}\boldsymbol{V}\right)
	\end{array}\label{eq:Phi_B}
\end{equation}
with $\boldsymbol{S}_{B}=\sum_{j=1}^{c}\left(\boldsymbol{M}_{j}-\boldsymbol{M}\right)\left(\boldsymbol{M}_{j}-\boldsymbol{M}\right)^{\top}\left(\in\mathbb{C}^{n_{1}n_{3}\times n_{1}n_{3}}\right)$.
Using the same steps for $\psi_{W}(\mathcal{V})$ , we obtain
\begin{equation}
	\widetilde{\psi}_{W}(\boldsymbol{V})=\mathrm{Trace}\left(\boldsymbol{V}^{\top}\boldsymbol{S}_{W}\boldsymbol{V}\right)\label{eq:Phi_W}
\end{equation}
with $\boldsymbol{S}_{W}=\sum_{j=1}^{c}\sum_{i\in C_{j}}\left(\boldsymbol{X}_{i}-\boldsymbol{M}_{j}\right)\left(\boldsymbol{X}_{i}-\boldsymbol{M}_{j}\right)^{\top}\left(\in\mathbb{C}^{n_{1}n_{3}\times n_{1}n_{3}}\right)$.\\
It can be noticed that $\boldsymbol{S}_{B}$ (or $\boldsymbol{S}_{W}$)
represents a block diagonal matrix where the $i^{th}$ block is the
the frontal slice $\widetilde{S}_{B}^{(i)}$ (or $\widetilde{S}_{W}^{(i)}$)
of the third-order tensor $\widetilde{\mathcal{S}}_{B}$ (or $\mathcal{\widetilde{S}}_{W}$). Then a new objective function equivalent to (\ref{eq:C-1}) is
defined in the transform domain by
\begin{equation}
	\boldsymbol{V}^{*}\in\mathbb{R}^{n_{1}n_{3}\times Kn_{3}}:\underset{\boldsymbol{V}}{max}\;\frac{\widetilde{\psi}_{B}(\boldsymbol{V})}{\widetilde{\psi}_{W}(\boldsymbol{V})}. \label{eq:C-2}
\end{equation}
As in the matrix case, (\ref{eq:C-2}) can be solved either by the Newton-Lanczos
algorithm (\ref{eq:LDA_3}) or by eigen-decomposition with regularization (\ref{eq:LDA_6}). The Newton-Lanczos algorithm involves iteratively the eigen-value decomposition of the matrix $\boldsymbol{S}(\rho)=\boldsymbol{S}_{B}-\rho\boldsymbol{S}_{W}$.
Since $\boldsymbol{S}(\rho)$ is a block diagonal matrix, i.e. $\boldsymbol{S}(\rho)=bdiag\left(\widetilde{S}^{(i)}\left(\rho\right)\right)$
where $\widetilde{S}\left(\rho\right)\in\mathbb{C}^{n_{1}\times n_{1}}$,
for $i=1,...,n_{3}$, eigen-value decomposition can be computed on
each block separately. Concerning the regularized eigen-decomposition problem (\ref{eq:LDA_3}), the inversion of the block diagonal matrix $\boldsymbol{S}_{W}(\gamma)=(\boldsymbol{S}_{W}+\gamma \boldsymbol{I})=bdiag(\widetilde{S}_{W}^{(i)}(\gamma))$ is also a block diagonal matrix where each block $\widetilde{S}_{W}^{(i)}$ is separately inverted. Since $\boldsymbol{S}_{W}(\gamma)$ and  $\boldsymbol{S}_{B}$ are square matrices identically partitioned into block diagonal form, the product $\boldsymbol{S}_{W}^{-1}(\gamma)\boldsymbol{S}_{B}$ also forms a diagonal block matrix identically partitioned. Eigen-decomposition can then be  computed on each block separately. Algorithms 2 and Algorithm 3 summarize the main steps for computing
$*_{L}$-TLDA either formulated as the trace ratio problem or the ratio trace problem.

\begin{algorithm}[!h] 	
	\caption{$*_{L}$TLDA - Trace ratio optimization}\label{algo2}
	\begin{algorithmic} 
		\STATE \textbf{Inputs}: $\mathcal{X}\in\mathbb{R}^{n_{1}\times K\times n_{3}}$ (input data:third-order tensor),\\
		\STATE \qquad \qquad$Y$ (labels: c classes) \\
		\STATE \qquad \qquad$K$ (reduced dimension) \\
		\STATE \textbf{Output}: $\mathcal{V}^{*}\in\mathbb{R}^{n_{1}\times K\times n_{3}}$ (projective tensor)\\
		\STATE \textbf{       }
		\STATE $\widetilde{\mathcal{X}}=L(\mathcal{X})$
		\FOR{$i=1,\ldots,n_{3}$} 
		\STATE $\widetilde{S}_{W}^{(i)},\widetilde{S}_{B}^{(i)}$ $\leftarrow$BuildScatters($\widetilde{X}^{(i)}$,$Y$)
		\STATE $(\widetilde{V}^{*})^{(i)}$ $\leftarrow$NewtonLanczos($\widetilde{S}_{W}^{(i)}$,$\widetilde{S}_{B}^{(i)}$,$K$) \qquad (see Algorithm 1)
		\ENDFOR \\
		$\mathcal{V}^{*}=L^{-1}(\widetilde{\mathcal{V}}^{*})$
	\end{algorithmic}
\end{algorithm}

\begin{algorithm}[!h] 	
	\caption{$*_{L}$TLDA - Ratio trace optimization}\label{algo3}
	\begin{algorithmic} 
		\STATE \textbf{Inputs}: $\mathcal{X}\in\mathbb{R}^{n_{1}\times K\times n_{3}}$ (input data:third-order tensor),\\
		\STATE \qquad \qquad$Y$ (labels: c classes) \\
		\STATE \qquad \qquad$\gamma>0$ (regularization parameter)\\
		\STATE \textbf{Output}: $\mathcal{V}^{*}\in\mathbb{R}^{n_{1}\times K\times n_{3}}$ (projective tensor)\\
		\STATE \textbf{       }
		\STATE $\widetilde{\mathcal{X}}=L(\mathcal{X})$
		\FOR{$i=1,\ldots,n_{3}$} 
		\STATE $\widetilde{S}_{W}^{(i)},\widetilde{S}_{B}^{(i)}$ $\leftarrow$ BuildScatters($\widetilde{X}^{(i)}$,$Y$)
		\STATE $K \leftarrow rank(S_B)$
		\STATE $(\widetilde{V}^{*})^{(i)} \leftarrow$ eigs$((\widetilde{S}_{W}^{(i)}+\gamma I)^{-1}\widetilde{S}_{B}^{(i)},$K,'lm')
		\ENDFOR \\
		$\mathcal{V}^{*}=L^{-1}(\widetilde{\mathcal{V}}^{*})$
	\end{algorithmic}
\end{algorithm}

\section{Experimental Results}

\subsection{Data sets} 
\begin{figure}[htp]
	\centering
	\subfloat[DIV (dim:64×64×2×5000) ]{\label{figur:1}\includegraphics[width=5cm,height=4cm]{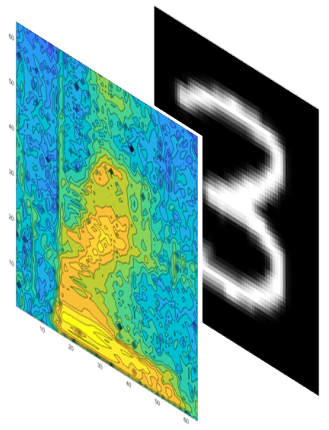}}
	\subfloat[WDCM (dim: 7×7×191×8032) ]{\label{figur:2}\includegraphics[width=5cm,height=4cm]{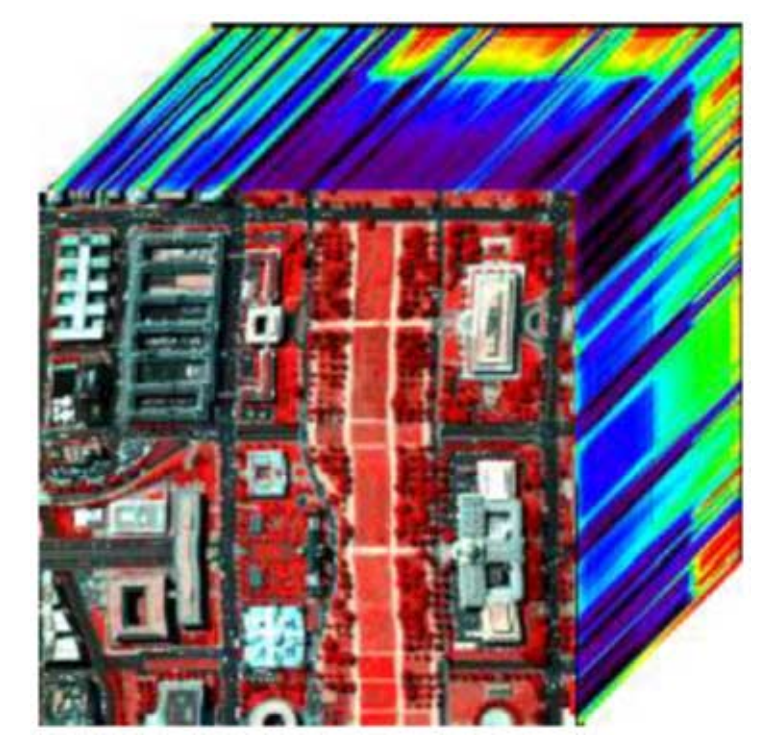}}
	\\
	\subfloat[AR (dim:25×25×3×2600) ]{\label{figur:3}\includegraphics[width=5cm,height=4cm]{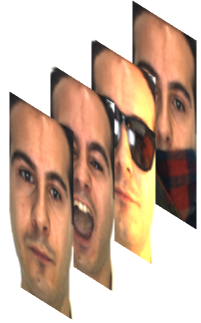}}
	\subfloat[GAIT (dim:32×32×10×4527) ]{\label{figur:4}\includegraphics[width=5cm,height=4cm]{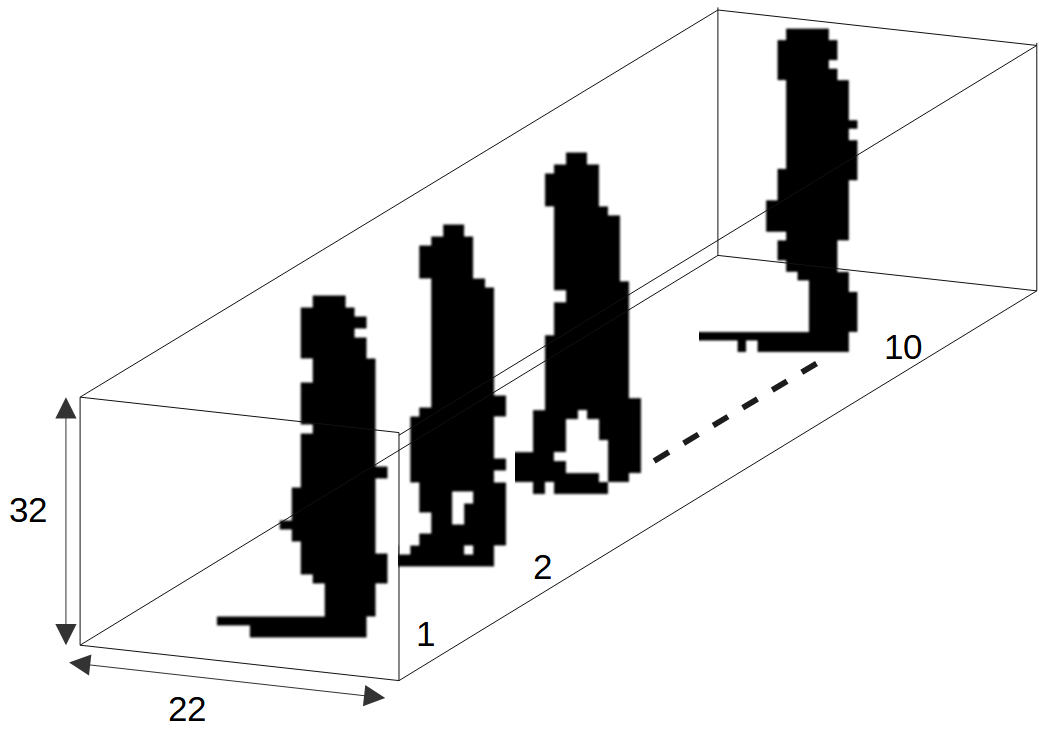}}
	
	\label{figur}\caption{Data sets}
\end{figure}
\begin{table}[htp]
	\caption{Gait data set : characteristics}
	\centering 
	\begin{tabular}{ccccc} \hline  \noalign{\vskip0.1cm} 
		Data sets &Gallery (GAR)& A (GAL) & B (GBR)& C (GBL)\tabularnewline \hline  \noalign{\vskip0.2cm} 
		
		nb of seq. &731&727&422&419 \tabularnewline
		nb of subj. &71&71&41&41 \tabularnewline \noalign{\vskip0.2cm} 
		Diff. from the gallery set& - &View&Shoe&Shoe-view \tabularnewline
		\noalign{\vskip0.2cm} \hline 
		
	\end{tabular} 
\end{table}

The experimental evaluation is based on five multidimensional data sets: The Digit Image Voice (DIV) data set, the Facial Recognition Technology (FERET) database , the AR Face (AR) database, the Washington DC Mall (WDCM) dataset  and the HumanID gait (GAIT) data set. Figure 5.1 illustrates four examples of the studied data sets. 

\textbf{The WDCM dataset} is a 191 band hyperspectral image of Washington DC Mall collected by the Hyperspectral Digital Imagery Collection Experiment (HYDICE) \cite{He2015}\footnote{ http://lesun.weebly.com/hyperspectral-data-set.html}. The whole image contains 1208×307 pixels. From this image we defined 4 classes: 'grass land', 'tree', 'roof' and 'road' which are hand-built and defined by 7×7 image blocs of pixels. We collected a total of 8032 blocs: 1894 for the 'tree', 1919 for the 'grass land', 2616 for the 'roof' and 1603 for the 'road'. A tensor representation is then built and defined by a 7×7×191×8032.

\textbf{The DIV dataset} build a tensor representation of the digits 0-9 obtained from two modalities: visual and audio and based on the MNIST \footnote{http://yann.lecun.com/exdb/mnist} and FSDD datasets \footnote{https://github.com/Jakobovski/free-spoken-digit-dataset}, respectively. The MNIST dataset contains 60000 training and 10000 test grayscale images of handwritten digits, with dimensions of 28x28 pixels. The FSDD dataset consists of 500, 8 kHz recordings of English pronunciations of the digits 0-9. These recordings are of varying durations, with a mean of approximately 0.5s. We preprocessed the recordings by converting them into 64×64 grayscale spectrograms. In order to harmonize the image sizes between modalities, the MNIST images are resized to 64×64 pixels. A tensor representation is then generated from 5000 samples combining the two modalities where each sample is randomly selected both in the resized MNIST dataset and the FSDD data sets. We obtain a 4th-order tensor of size 64×64×2×5000.

\textbf{The FERET database} is a standard facial image collection including 14126 images from 1199 individuals with different view points \cite{Phillips2000}\footnote{https://old.datahub.io/dataset/feret-database}. In our experiment, we select a subset composed of 80 subjects where each of them having at least 10 images, resulting in 1145 images. We generate a 3th-order tensor of size 32×32×1145.

\textbf{The AR dataset} contains over 4,000 color images corresponding to 126 people's faces (70 men and 56 women) \cite{Mart1998} \footnote{http://cbcsl.ece.ohio-state.edu/ARdatabaseNew.html}.  Images feature frontal view faces with different facial expressions, illumination conditions, and occlusions (sun glasses and scarf). The pictures were taken under strictly controlled conditions. No restrictions on wear (clothes, glasses, etc.), make-up, hair style, etc. were imposed to participants. A subset of 100 subjects have been considered corresponding to a total of 2600 images of size 25×25. A tensor of size 25×25×3×2600 was generated.

\textbf{The GAIT dataset} is build from the USF HumanID “Gait Challenge” data sets version 1.7 \cite{Sarkar2005}\footnote{http://www.eng.usf.edu/cvprg/Gait-Data.html}. This data set is composed of 452 sequences from 74 subjects walking in elliptical paths in front of the camera. For each subject, there are three covariates: viewpoint (left/right), shoe type (two different types) and surface type (grass/concrete). In our experiments, we consider only the sequences corresponding to the \textit{grass} type surface defining thus the “gallery” set. This dataset contains 731 sequences from 71 subjects (persons) and each subject has an average of roughly 10 samples available under the form of binary silhouette images of size 32×32 (see Figure 5.1). Thus we define a 4th-order training tensor of size 32×32×10×731. The test set is  based on three probe sets named A, B and C as detailed in  Table 5.1.  More precisely, the image acquisition conditions for the “gallery set” and each probe set are summarized in brackets after the data name in Table 5.1, where G, A, B, L, and R stand for grass surface, shoe type A, shoe type B, left view, and right view, respectively. There is no redundancy between the gallery set and each probe set, i.e. there are no common subjects and sequences between them. 

\begin{table} 
	\caption{Performances on DIV dataset (9 categories).}
	\centering 
	\begin{tabular}{lllll} \hline  \noalign{\vskip0.1cm} 
		Objective & Methods & ACC ($\%$) & Times (sec)& DIM\tabularnewline \hline  \noalign{\vskip0.01cm} 
		\multirow{6}{*}{trace ratio} 
		\tabularnewline  & Fisherfaces &83.08($\pm$0.026)&0.33 ($\pm$0.013)&9
		\tabularnewline  & *c-TDA &\cellcolor{gray!15}90.58 ($\pm$0.02)&0.45 ($\pm$0.018)&18
		\tabularnewline  & *t-TDA &\cellcolor{gray!5}89.75 ($\pm$0.02)&0.46 ($\pm$0.026)&18
		\tabularnewline  & DATER &87.98 ($\pm$0.02)&0.33 ($\pm$0.016)&928
		\tabularnewline  & CMDA &89.43($\pm$0.023)&0.42 ($\pm$0.058)&493
		\tabularnewline  & UMDA &74.95 ($\pm$0.04)&11.44 ($\pm$0.08)&28
		\tabularnewline \noalign{\vskip0.05cm} 
		\multirow{5}{*}{ratio trace}
		\tabularnewline  & *c-TDA&89.33($\pm$0.017)&1.79 ($\pm$1.9)&823
		\tabularnewline  & *t-TDA &88.91 ($\pm$0.02)&1.50($\pm$0.39)&815
		\tabularnewline  & DGTDA&89.13($\pm$0.017)&0.34 ($\pm$0.097)&968
		\tabularnewline  & HODA&87.46($\pm$0.022)&0.64 ($\pm$0.07)&956
		\tabularnewline \noalign{\vskip0.1cm} \hline 
		
	\end{tabular} 
\end{table}
\begin{table} 
	\caption{Performances on WDCM dataset (4 categories)}
	\centering 
	\begin{tabular}{lllll} \hline  \noalign{\vskip0.1cm} 
		Objective & Methods & ACC ($\%$) & Times (sec)& DIM\tabularnewline \hline  \noalign{\vskip0.01cm} 
		\multirow{6}{*}{trace ratio} 
		\tabularnewline  &Fisherfaces &91.6($\pm$0.019)&0.478 ($\pm$0.05)&3
		\tabularnewline  & *c-TDA &\cellcolor{gray!5}95.05 ($\pm$0.01)&0.45 ($\pm$0.010)&384
		\tabularnewline  & *t-TDA &\cellcolor{gray!15}95.51 ($\pm$0.01)&0.61 ($\pm$0.01)&564
		\tabularnewline  & DATER &94.21($\pm$0.016)&2.28 ($\pm$0.017)&150
		\tabularnewline  & CMDA &96.65($\pm$0.01)&2.11 ($\pm$0.054)&305
		\tabularnewline  & UMDA &91.95 ($\pm$0.025)&118.08 ($\pm$3.52)&28
		\tabularnewline \noalign{\vskip0.05cm}  
		\multirow{5}{*}{ratio trace} 
		\tabularnewline  & *c-TDA&93.78($\pm$0.04)&0.75 ($\pm$0.08)&964
		\tabularnewline  & *t-TDA &94.41 ($\pm$0.02)&1.01($\pm$0.15)&1000
		\tabularnewline  & DGTDA&94.35($\pm$0.017)&2.60 ($\pm$0.09)&350
		\tabularnewline  & HODA&94.45($\pm$0.044)&4.42 ($\pm$0.77)&783
		\tabularnewline \noalign{\vskip0.1cm} \hline 
		
	\end{tabular} 
\end{table}
\begin{table} 
	\caption{Performances on FERET dataset (80 categories).}
	\centering 
	\begin{tabular}{lllll} \hline  \noalign{\vskip0.1cm} 
		Objective & Methods & ACC ($\%$) & Times (sec)& DIM\tabularnewline \hline  \\ 
		\multirow{6}{*}{trace ratio} 
		\tabularnewline  &Fisherfaces &\cellcolor{gray!15}88.71($\pm$0.02)&0.61 ($\pm$0.02)&79
		\tabularnewline  & *c-TDA &\cellcolor{gray!5}87.93 ($\pm$0.02)&1.51 ($\pm$0.02)&68
		\tabularnewline  & *t-TDA &87.78 ($\pm$0.02)&1.54 ($\pm$0.016)&68
		\tabularnewline  & DATER &81.23 ($\pm$0.03)&1.11 ($\pm$0.19)&45
		\tabularnewline  & CMDA &75.46($\pm$0.02)&1.32 ($\pm$0.02)&559
		\tabularnewline  & UMDA &78.2 ($\pm$0.03)&5.85 ($\pm$0.27)&28
		\tabularnewline  \noalign{\vskip0.05cm}   
		\multirow{5}{*}{ratio trace} 
		\tabularnewline  & *c-TDA&63.75($\pm$0.03)&6.08 ($\pm$1.24)&675
		\tabularnewline  & *t-TDA &64.01 ($\pm$0.03)&5.50($\pm$1.44)&527
		\tabularnewline  & DGTDA&64.46($\pm$0.04)&0.487 ($\pm$0.02)&607
		\tabularnewline  & HODA&67.96($\pm$0.02)&1.38($\pm$0.61)&29
		\tabularnewline \hline 
		
	\end{tabular} 
\end{table}
\begin{table} 
	\caption{Performances on AR dataset (100 categories).}
	\centering 
	\begin{tabular}{lllll} \hline  \noalign{\vskip0.1cm} 
		Objective & Methods & ACC ($\%$) & Times (sec)& DIM\tabularnewline \hline  \\ 
		\multirow{6}{*}{trace ratio} 
		\tabularnewline  &Fisherfaces &92.25($\pm$0.014)&0.51 ($\pm$0.021)&99
		\tabularnewline  & *c-TDA &\cellcolor{gray!5}94.36 ($\pm$0.02)&1.49 ($\pm$0.04)&294
		\tabularnewline  & *t-TDA &\cellcolor{gray!15}95.58 ($\pm$0.01)&2.92 ($\pm$0.08)&240
		\tabularnewline  & DATER &68.9 ($\pm$0.03)&4.44 ($\pm$0.16)&101
		\tabularnewline  & CMDA &72.23($\pm$0.03)&5.13 ($\pm$0.17)&164
		\tabularnewline  & UMDA &68.2 ($\pm$0.035)&13.46 ($\pm$0.73)&23
		\tabularnewline  \noalign{\vskip0.05cm}    
		\multirow{5}{*}{ratio trace} 
		\tabularnewline  & *c-TDA&33.86($\pm$0.03)&2.31 ($\pm$0.05)&956
		\tabularnewline  & *t-TDA &34.48 ($\pm$0.03)&4.21($\pm$0.36)&438
		\tabularnewline  & DGTDA&33.83($\pm$0.02)&1.47 ($\pm$0.077)&892
		\tabularnewline  & HODA&34.88($\pm$0.04)&1.54 ($\pm$0.076)&984
		\tabularnewline \hline 
		
	\end{tabular} 
\end{table}

%

\subsection{Competitors} 
We compare our approach with five supervised learning algorithms:\\ PCA+LDA (Fisherfaces) \cite{Belhumeur1997}, Discriminant analysis with tensor representation (DATER)  \cite{Yan,Visani2005}, Constrained Multilinear Discriminant Analysis (CMDA) \cite{li}, Direct General Tensor Discriminant Analysis (DGTDA) \cite{li}, Higher Order Discriminant Analysis (HODA)  \cite{Phan2010}, Uncorrelated Multilinear Discriminant Analysis with regularization (UMLDA)  \cite{Lu2009}.
Our approach will be tested with two tensor-tensor products: t-product and the c-product. When using the t-product, our approach will be refered as $_{*t}$-TDA and $_{*c}$-TDA for the c-product.

\noindent Fisherfaces's implementation is based on the ratio trace criteria and uses a vector to vector projection. In order to avoid the singularity problem of the within-class scatter matrix, PCA is beforehand computed reducing thus the dimension of the feature space. To set the output dimension, a classical heuristic is to retain the k eigenvectors that capture a certain percentage of the total variance. In all the experiments, we will consider at least 95$\%$ of the total variance.  DATER, CMDA, DGTDA, HODA and UMLDA are multidimensional variants of LDA. The first four methods use a tensor to tensor projection while the last one uses a tensor to vector projection. DATER, CMDA and UMLDA are formulated as a trace ratio problem while DGTDA and HODA solve the ratio trace problem. \\ 
When the optimization problem is formulated as the trace ratio problem, a regularization parameter is used in order to avoid the SSS problem (singularity problem of scatter matrices). The regularization parameter is selected by a k-fold cross validation step. As for the dimension of the projective subspace, it is given by the rank of the between-scatter matrix $S_B$. When the optimization problem is formulated as the ratio trace problem, the dimension of the projective subspace is determined  by k-fold cross validation. The tests are based on 30 repetitions of the experiments and the average accuracy is used as classification performance. \\

\noindent Tables 5.2--5.5 summarize the performances of the different methods when applied to DIV, WDCM, FERET and AR data sets, respectively. These tables show the values of the average accuracy (ACC) recorded by the methods (third column), the average training times (fourth column) and the maximal output dimensions (fifth column). First, from a general view, we observe that the proposed $*_{L}$-TLDA provide similar results whether it is based on the c-product or the t-product. Secondly, when formulated as the trace ratio problem, the proposed method clearly records better classification results with the other tensor decomposition methods based on the same objective such as DATER, CMDA and UMDA. It also outperforms the version based on the ratio trace objective and very clearly when the number of categories increases as in the AR data base (Table 5.5) and the FERET data base (Table 5.4). This result is also valid for the other competitors based on the ratio trace objective such as DGTDA and HODA.\\
Moreover, we observe that $*_{L}$-TLDA (trace ratio criterion) records competitive training times with the other studied methods. By using a tensor to vector projection strategy, UMLDA shows the highest complexity making this method very time consuming.

In a second experiment, we study the performance of the  proposed algorithm on the Gait sequence.  The identification performance is measured by the Cumulative Match Characteristic (CMC) as defined in \cite{Boulgouris2005}  which plots identification rates within a given rank k. More precisely, rank k results report the percentage of probe subjects whose the true match in the gallery set was in the top k matches.  The rank 1 and the rank 5 gait recognition results using the modified angle distance (MAD) \cite{Lu2008} are presented in Table 5.6. As previously, $*_{L}$-TLDA formulated with the trace ratio criterion shows the best recognition rates on all the probe sets. When the ratio trace criterion is optimized, the results are markedly lower. This result confirms those of Tables 5.2--5.5. We can notice that Fisherfaces records the lowest recognition rates making it clear that a matricial treatment of this kind of data set is not well suited.
\begin{table} 
	\caption{Performances on the Gait sequence (72 categories).}
	\centering 
	\begin{tabular}{cccccccc} \hline  \noalign{\vskip0.1cm} 
		Objective&Methods& \multicolumn{3}{c}{rank 1} & \multicolumn{3}{c}{rank 5} \tabularnewline \noalign{\vskip0.05cm} \hline
		
		& & A & B & C & A & B &C \tabularnewline  \noalign{\vskip0.05cm} \cline{3-8}  
		\multirow{6}{*}{TRO} 
		\tabularnewline  &Fisherfaces&43.6&43.9&31.7&71.8&63.4&51.20
		\tabularnewline  & *c-TDA &\cellcolor{gray!5}90.1&\cellcolor{gray!15}84.9&\cellcolor{gray!15}65.8&\cellcolor{gray!15}100&\cellcolor{gray!15}92.7&\cellcolor{gray!5}85.4
		\tabularnewline  & *t-TDA &\cellcolor{gray!15}94.4&\cellcolor{gray!5}80.5&\cellcolor{gray!15}65.8&\cellcolor{gray!15}100&\cellcolor{gray!15}92.7&\cellcolor{gray!15}87.8
		\tabularnewline  & DATER &69.0&70.7&48.8&93.0&82.9&68.3
		\tabularnewline  & CMDA &80.2&73.7&53.6&95.8&80.5&70.7
		\tabularnewline  & UMDA &77.5&68.3&41.5&94.4&80.5&73.2
		\tabularnewline  \noalign{\vskip0.05cm}    
		\multirow{5}{*}{RTO} 
		\tabularnewline  & *c-TDA&65.2&69.3&46.1&92.4&80.5&67.2
		\tabularnewline  & *t-TDA &67.8&72.3&47.2&94.2&81.2&69.5
		\tabularnewline  & DGTDA&63.4&73.2&44.0&90.1&82.9&68.3
		\tabularnewline  & HODA&70.4&70.7&48.8&92.9&80.5&68.3
		\tabularnewline \hline 
		
	\end{tabular} 
\end{table}

\subsection{Conclusion}
In this work, we proposed a new Tensor Linear Discriminant Analysis based on the the concept of transform domain as defined in \cite{Kernfeld2015}. Considering any fixed, invertible linear transformation $L$, our $_{*L}$TDA procedure is based on the computation of a new  $_{*L}$-family product such as the t-product or the c-product. In this context, we showed that the solution of our MDA can be obtained in a fully tensor form instead of a sequence of projective matrices, solutions of existing MDA methods. Another key aspect is that the obtained solution is the result of independent optimization problems easier to solve and more robust compared to existing MDA methods that  are based on alternating optimization heuristics. The experimental evaluation based on these two products show similar classification performances with a slight advantage of the c-product in terms of training time. The experimental evaluations show that the choice of the optimization criterion, i.e. the trace ratio or the ratio trace, influences significantly the classification performances of our method. The conclusions of our experimental evaluation $_{*L}$TDA, based on the trace ratio criterion performs very well and outperforms most of the existing MDA methods. 

Several issues remain to be investigated. First, the proposed MDA is based on the   building of three-order tensors and its extension to higher-order tensors could be the subject of future work.
Second, the concept of rank being clearly defined in traditional Linear Discriminant Analysis, it will be interesting to address this issue shortly in the framework of $_{*L}$TDA in order to bound the dimensionality of the solution. 

%
%
%
%
%

\end{document}